\documentclass[12pt,savetrees]{amsart}
\usepackage{amssymb, hyperref}

\newcommand{\Bgp}{{\Z^\N}}

\long\def\forget#1\forgotten{}
\newcommand{\issuenumber}{34}
\newcommand{\issuemonth}{August}
\newcommand{\issueyear}{2012}

\newcommand{\ed}{
\newpage

\section{Unsolved problems from earlier issues}

\begin{issue}Is $\binom{\Omega}{\Gamma}=\binom{\Omega}{\Tau}$?\end{issue}%\stepcounter{issue}
\begin{issue}Is $\ufin(\cO,\Omega)=\sfin(\Gamma,\Omega)$?And if not, does $\ufin(\cO,\Gamma)$ imply
$\sfin(\Gamma,\Omega)$?\end{issue}%\stepcounter{issue}
\stepcounter{issue}\begin{issue}Does $\sone(\Omega,\Tau)$ imply $\ufin(\Gamma,\Gamma)$?\end{issue}
\begin{issue}Is $\fp=\fp^*$? (See the definition of $\fp^*$ in that issue.)\end{issue}
\begin{issue}Does there exist (in ZFC) an uncountable set satisfying $\sfin(\cB,\cB)$?\end{issue}
\stepcounter{issue}
\begin{issue}Does $X \nin \NON(\cM)$ and $Y\nin\mathsf{D}$ imply that $X\cup Y\nin \COF(\cM)$?\end{issue}
\begin{issue}[CH]Is $\split(\Lambda,\Lambda)$ preserved under finite unions?\end{issue}
\begin{issue}Is $\cov(\cM)=\fo$? (See the definition of $\fo$ in that issue.)\end{issue}
\stepcounter{issue}
\begin{issue}Could there be a Baire metric space $M$ of weight $\aleph_1$ and a partition
$\mathcal{U}$ of $M$ into $\aleph_1$ meager sets where for each ${\mathcal U}'\subset\mathcal U$,
$\bigcup {\mathcal U}'$ has the Baire property in $M$?\end{issue}
\stepcounter{issue} %% no problem in Issue 13
\begin{issue}Does there exist (in ZFC) a set of reals $X$ of cardinality $\fd$ such that all
finite powers of $X$ have Menger's property $\sfin(\cO,\cO)$?\end{issue}
\begin{issue}Can a Borel non-$\sigma$-compact group be generated by a Hurewicz subspace?\end{issue}
\begin{issue}[MA]Is there $X\sbst\bbR$ of cardinality continuum, satisfying $\sone(\BO,\BG)$?\end{issue}
\begin{issue}[CH]Is there a totally imperfect $X$ satisfying $\ufin(\cO,\Gamma)$
that can be mapped continuously onto $\Cantor$?\end{issue}
\begin{issue}[CH]Is there a Hurewicz $X$ such that $X^2$ is Menger but not Hurewicz?\end{issue}
\begin{issue}Does the Pytkeev property of $C_p(X)$ imply that $X$ has Menger's property?\end{issue}
\begin{issue}Does every hereditarily Hurewicz space satisfy $\sone(\BG,\BG)$?\end{issue}
\begin{issue}[CH]Is there a Rothberger-bounded $G\le\Bgp$ such that $G^2$ is not Menger-bounded?\end{issue}
\begin{issue}Let $\cW$ be the van der Waerden ideal. Are $\cW$-ultrafilters closed under products?\end{issue}
\begin{issue}Is the $\delta$-property equivalent to the $\gamma$-property $\binom{\Omega}{\Gamma}$?\end{issue}
\stepcounter{issue}\stepcounter{issue}
\general\end{document}}

\newcommand{\Cantor}{{\{0,1\}^\N}}

\newcommand{\fd}{\mathfrak{d}}
\newcommand{\fp}{\mathfrak{p}}
\newcommand{\NON}{{\mathsf   {NON}}}\newcommand{\COF}{{\mathsf   {COF}}}

\newcommand{\cM}{\mathcal{M}}
\newcommand{\cov}{\mathsf{cov}}

\newcommand{\bbR}{\mathbb{R}}
\newcommand{\fo}{\mathfrak{od}}

\newcommand{\w}{\omega}\newcommand{\ft}{\mathfrak{t}}
\renewcommand{\split}{\mathsf{Split}}\newcommand{\bq}{\begin{quote}}\newcommand{\eq}{\end{quote}}
\newcommand{\cO}{\mathcal{O}}\newcommand{\cB}{\mathcal{B}}\newcommand{\BG}{\cB_\Gamma}
\newcommand{\BO}{\cB_\Omega}

\newcommand{\sone}{\mathsf{S}_1}\newcommand{\sfin}{\mathsf{S}_\mathrm{fin}}
\newcommand{\ufin}{\mathsf{U}_\mathrm{fin}} 
\newcommand{\nin}{\not\in}
\newcommand{\cW}{\mathcal{W}}

\newcommand{\N}{\mathbb{N}}\newcommand{\Z}{\mathbb{Z}}
\newcommand{\sbst}{\subseteq}
\newcommand{\by}[2]{\par\hfill\emph{#1}, #2}\newcommand{\Tau}{\mathrm{T}}
\newcommand{\CE}{\textsc{CE}}
\newtheorem{thm}{Theorem}[section]\newcommand{\bthm}{\begin{thm}} \newcommand{\ethm}{\end{thm}}
\newtheorem{prop}[thm]{Proposition}\newcommand{\bprp}{\begin{prop}} \newcommand{\eprp}{\end{prop}}
\newtheorem{fact}[thm]{Fact}\newcommand{\bfct}{\begin{fact}} \newcommand{\efct}{\end{fact}}
\newtheorem{prob}[thm]{Problem}\newcommand{\bprb}{\begin{prob}} \newcommand{\eprb}{\end{prob}}
\newtheorem{lem}[thm]{Lemma}\newcommand{\blem}{\begin{lem}} \newcommand{\elem}{\end{lem}}
\newtheorem{claim}[thm]{Claim}\newcommand{\bclm}{\begin{claim}} \newcommand{\eclm}{\end{claim}}
\newtheorem{cor}[thm]{Corollary}\newcommand{\bcor}{\begin{cor}} \newcommand{\ecor}{\end{cor}}
\newtheorem{conj}[thm]{Conjecture}\newcommand{\bcnj}{\begin{conj}} \newcommand{\ecnj}{\end{conj}}
\theoremstyle{definition}\newtheorem{defn}[thm]{Definition}\newcommand{\bdfn}{\begin{defn}} \newcommand{\edfn}{\end{defn}}
\theoremstyle{remark}\newtheorem{rem}[thm]{Remark}\newcommand{\brem}{\begin{rem}} \newcommand{\erem}{\end{rem}}
\newtheorem{cnv}[thm]{Convention}\newcommand{\bcnv}{\begin{cnv}} \newcommand{\ecnv}{\end{cnv}}
\newtheorem{exam}[thm]{Example}\newcommand{\bexm}{\begin{exam}} \newcommand{\eexm}{\end{exam}}
\newtheorem{issue}{Issue}\newcommand{\bpf}{\begin{proof}} \newcommand{\epf}{\end{proof}}
\newcommand{\be}{\begin{enumerate}}\newcommand{\ee}{\end{enumerate}}\newcommand{\bi}{\begin{itemize}}
\newcommand{\ei}{\end{itemize}}
\newcommand{\general}{\small\vfill\par\noindent\hrulefill\par
\noindent\textbf{Previous issues.} The previous issues of this
bulletin are available online at\\
\url{http://front.math.ucdavis.edu/search?\&t=\%22SPM+Bulletin\%22}
\\[0.1cm]
\textbf{Contributions.} Announcements, discussions, and open problems should be emailed
to \texttt{tsaban@math.biu.ac.il}\\[0.1cm]
\textbf{Subscription.}
To receive this bulletin (free) to your e-mailbox, e-mail us.
}

\newcommand{\arXivl}[4]{\subsection{#2}{#4}\par\hfill{\arx{#1}}\par\hfill\emph{#3}}
\newcommand{\arXiv}[3]{\subsection{#2}\mbox{}\par\hfill{\arx{#1}}\par\hfill\emph{#3}}

\newcommand{\nAMSPaper}[4]{\subsection{#2}{#4}\par\hfill{\texttt{#1}}\par\hfill\emph{#3}}
\newcommand{\AMS}[3]{\subsection{#1}\mbox{}\par\hfill{\texttt{#3}}\par\hfill\emph{#2}}

\newcommand{\arx}[1]{\url{http://arxiv.org/abs/#1}}

\title[$\mathcal{SPM}$ Bulletin \textbf{\issuenumber} (\issuemonth{} \issueyear)]{%
$\mathcal{SPM}$ Bulletin\\[0.5cm]
Issue number \issuenumber: \issuemonth{} \issueyear{} \CE{}}

\begin{document}
\maketitle

%\tableofcontents

\section{Editor's note}

Dear Friends,

\medskip

1. Marion Scheepers, of the Department of Mathematics,
Boise State University,
was awarded the honor and title of Distinguished Professor.
This award is, among other things, in recognition of Prof.\ Scheepers' founding the field
of Selection Principles and of his seminal contributions to this field.
I use this opportunity to greet Prof. Scheepers for this well-deserved recognition, and wish
him many more happy years of fruitful research.

\medskip

2. Bemaer presentations for most of the lectures delivered at that
\emph{Fourth Workshop on Coverings, Selections and Games in Topology}
(aka SPMC2012),
are available at the Workshop's webpage:
\url{http://u.cs.biu.ac.il/~tsaban/spmc12/}
under \emph{Lectures}.

\medskip

3. A special Issue of \emph{Topology and its Applications} will be dedicated to SPM and related topics,
following the SPMC2012 conference. Submissions of high quality research papers is welcome.
Submit your paper to \texttt{TopApplIssue@gmail.com}

Submissions will be refereed according to the usual high standards of the journal, and the final
acceptance or rejection decision will be made by the journal's chief editors.

\textbf{Submission deadline:} October 31, 2012. Later submissions may be considered in
exceptional cases. Email the mentioned address in case of a need for deadline extension.

\medskip

4. Last, but not least: Our \emph{long announcements} section is concluded by one where
Malliaris and Shelah announce a \emph{solution to the minimal tower problem}
Surprisingly, their solution is that $\fp=\ft$, in ZFC! This solution is expected to have
deep implications within SPM and related areas. 
Details about the minimal tower problem are available at Issue 5 of this bulletin:
\url{http://arxiv.org/pdf/math/0305367.pdf}
Greetings to Malliaris and Shelah for their
breakthrough.

\medskip

With best regards,

\by{Boaz Tsaban}{tsaban@math.biu.ac.il}

\hfill \texttt{http://www.cs.biu.ac.il/\~{}tsaban}

\section{Long announcements}

\arXivl{1202.2056}
{Universally Kuratowski--Ulam space and open-open games}
{Andrzej Kucharski}
{We examine the class of spaces in which the second player has winning
strategy in the open-open game. We shown that this space is not an universally
Kuratowski-Ulam. We also show that the games G and G7 introduced by Daniels,
Kunen, Zhou are not equivalent.}

\nAMSPaper{http://www.ams.org/journal-getitem?pii=S0002-9947-2012-05402-X}
{Maximal functions and the additivity of various families of null sets}
{Juris Steprans}
{It is shown to be consistent with set theory that every set of reals of size $\aleph_1$
is null yet there are  planes in Euclidean 3-space whose union is not null.
Similar results are obtained for circles in the plane as well as other geometric objects.
The proof relies on results from harmonic analysis about the boundedness of certain maximal
operators and a measure-theoretic pigeonhole principle.}

\arXivl{1205.6824}
{Additivity of the Gerlits--Nagy property and concentrated sets}
{Boaz Tsaban and Lyubomyr Zdomskyy}
{We settle all problems posed by Scheepers, in his tribute paper to Gerlits,
concerning the additivity of the Gerlits--Nagy property and related additivity
numbers. We apply these results to compute the minimal number of concentrated
sets of reals (in the sense of Besicovitch) whose union, when multiplied with a
Gerlits--Nagy space, need not have Rothberger's property. We apply these
methods to construct a large family of spaces, whose product with every
Hurewicz space has Menger's property.}

\arXivl{1207.0827}
{On $\infty$-convex sets in spaces of scatteredly continuous functions}
{Taras Banakh, Bogdan Bokalo, and Nadiya Kolos}
{Given a topological space $X$, we study the structure of $\infty$-convex
subsets in the space $SC_p(X)$ of scatteredly continuous functions on $X$. Our
main result says that for a topological space $X$ with countable strong fan
tightness, each potentially bounded $\infty$-convex subset $F\subset SC_p(X)$
is weakly discontinuous in the sense that each non-empty subset $A\subset X$
contains an open dense subset $U\subset A$ such that each function $f|U$, $f\in
F$, is continuous. This implies that $F$ has network weight $nw(F)\le nw(X)$.}

\arXivl{1207.4025}
{Diagonalizations of dense families}
{Maddalena Bonanzinga, Filippo Cammaroto, Bruno Antonio Pansera, Boaz Tsaban}
{We develop a unified framework for the study of properties involving
diagonalizations of dense families in topological spaces. We provide complete
classification of these properties. Our classification draws upon a large
number of methods and constructions scattered in the literature, and on some
novel results concerning the classical properties.

\emph{Comment.} The field is very active these days, and consequently we may have
  missed some references. We would appreciate any feedback, in particular
  concerning relevant references.
}

\arXivl{1208.1182}
{On open-open games of uncountable length}
{Andrzej Kucharski}
{The aim of this note is to investigate the open-open game of uncountable
length. We introduce a cardinal number $\mu(X)$, which says how long the Player
I has to play to ensure a victory. It is proved that
$su(X)\leq\mu(X)\leq su(X)^+$. We also introduce the class $\mathcal C_\kappa$
of topological spaces that can be represented as the inverse limit of
$\kappa$-complete system $\{X_\sigma,\pi^\sigma_\rho,\Sigma\}$ with
$\w(X_\sigma)\leq\kappa$ and skeletal bonding maps. It is shown that product of
spaces which belong to $\mathcal C_\kappa$ also belongs to this class and
$\mu(X)\leq\kappa$ whenever $X\in\mathcal C_\kappa$ .
\par
Published in: International Journal of Mathematics and Mathematical Sciences,
  vol.\ 2012, Article ID 208693, 2012.}

\arXivl{1208.4823}
{Topologically invariant $\sigma$-ideals on Euclidean spaces}
{Taras Banakh, Micha\l{} Morayne, Robert Ra\l{}owski, Szymon \.Zeberski}
{We study and classify topologically invariant $\sigma$-ideals with an analytic
base on Euclidean spaces and evaluate the cardinal characteristics of such
ideals.}

\arXivl{1208.5424}
{Cofinality spectrum theorems in model theory, set theory and general
  topology}
{M. Malliaris and S. Shelah}
{We connect and solve two longstanding open problems in quite different areas:
the model-theoretic question of whether $SOP_2$ is maximal in Keisler's order,
and the question from set theory/general topology of whether $\mathfrak{p} =
\mathfrak{t}$, the oldest problem on cardinal invariants of the continuum. We
do so by showing these problems can be translated into instances of a more
fundamental problem which we state and solve completely, using model-theoretic
methods.}

\section{Short announcements}\label{RA}

\arXiv{1202.4018}
{Subspaces of monotonically normal compacta}
{Ahmad Farhat}

\AMS{The topological Baumgartner-Hajnal theorem}
{Rene Schipperus}
{http://www.ams.org/journal-getitem?pii=S0002-9947-2012-04990-7}

\AMS{Forcing, games and families of closed sets}
{Marcin Sabok}
{http://www.ams.org/journal-getitem?pii=S0002-9947-2012-05404-3}

\arXiv{1204.0918}
{Reflexivity in precompact groups and extensions}
{Monteserrat Bruguera, and Jorge Galindo, and Constancio Hern\'andez, and Mikhail Tkachenko}

\arXiv{1204.2438}
{On $\sigma$-convex subsets in spaces of scatteredly continuous functions}
{Taras Banakh, Bogdan Bokalo, Nadiya Kolos}

\arXiv{1204.3108}
{Filter convergence in $\beta\omega$}
{Jonathan L. Verner}

\AMS{Covering an uncountable square by countably many continuous functions}
{Wieslaw Kubis; Benjamin Vejnar}
{http://www.ams.org/journal-getitem?pii=S0002-9939-2012-11292-4}

\arXiv{1205.3391}
{Kuratowski operations revisited}
{Szymon Plewik and Marta Walczy\'nska}

\arXiv{1205.5909}{A new class of Ramsey-classification theorems and their applications in
 the Tukey theory of ultrafilters}{Natasha Dobrinen and Stevo Todorcevic}

\arXiv{1206.0795}
{On the structure of finite level and $\omega$-decomposable Borel functions}
{Luca Motto Ros}

\arXiv{1206.0722}
{Notes on the od-Lindel\"of property}
{Mathieu Baillif}

\AMS{Infinite dimensional perfect set theorems}
{Tamas Matrai}
{http://www.ams.org/journal-getitem?pii=S0002-9947-2012-05468-7}

\arXiv{1206.5581}{Large sets in the sense of the Ellentuck topology does not admit the
 Kuratowski's partition}{Ryszard Frankiewicz S{\l}awomir Szczepaniak}

\AMS{On Borel sets belonging to every invariant ccc $\sigma$-ideal on
$2^{\mathbb{N}}$}
{Piotr Zakrzewski}
{http://www.ams.org/journal-getitem?pii=S0002-9939-2012-11384-X}

\ed